\documentclass[11pt]{amsart}
\usepackage[T1]{fontenc}
\usepackage{amssymb,mathtools,booktabs}
\usepackage{mathrsfs}
\usepackage[margin=1.15in,top=0.79in,bottom=0.7in]{geometry}
\usepackage[all]{xy}
\usepackage{xcolor}

\newtheorem{maintheorem}{Theorem}
\newtheorem{computation}{Computation}
\newtheorem{theorem}{Theorem}[section]

\newtheorem{lemma}[theorem]{Lemma}
\newtheorem{corollary}[theorem]{Corollary}

\theoremstyle{remark}
\newtheorem{remark}[theorem]{Remark}

\numberwithin{equation}{section}

\DeclareMathOperator{\NS}{NS}
\DeclareMathOperator{\Pic}{Pic}
\DeclareMathOperator{\Br}{Br}
\DeclareMathOperator{\Gal}{Gal}
\DeclareMathOperator{\rk}{rk}
\DeclareMathOperator{\coker}{coker}
\DeclareMathOperator{\Res}{Res}

\newcommand{\ZZ}{\mathbb Z}
\newcommand{\QQ}{\mathbb Q}

\newcommand{\CC}{\mathbb C}
\newcommand{\FF}{\mathbb F}

\newcommand{\GG}{\mathbb G}
\newcommand{\Gm}{\GG_m}

\newcommand{\Qbar}{\overline{\QQ}}

\newcommand{\et}{\mathrm{\acute et}}
\newcommand{\tensor}{\otimes}
\newcommand{\isom}{\cong}

\newcommand{\linedef}[1]{\textsl{#1}}

\title[A cyclotomic Artin--Tate formula for K3 surfaces]
{A cyclotomic Artin--Tate formula and\\ sphere packing for K3 surfaces}

\author{Asher Auel}
\address{Department of Mathematics, Dartmouth College, Hanover, New Hampshire}
\email{asher.auel@dartmouth.edu}

\author{Sair Shaikh}
\email{sairshaikh1@gmail.com}


\usepackage{hyperref}
\hypersetup{
  pdftitle={A cyclotomic Artin--Tate formula for K3 surfaces},
  pdfauthor={Asher Auel and Sair Shaikh},
  colorlinks=true,
  linkcolor=blue,
  anchorcolor=blue,
  citecolor=blue
}
\begin{document}

\begin{abstract}
We give new necessary conditions on the Weil polynomial of a K3 surface over a finite field using techniques from the theory of sphere packing. By comparing the Artin--Tate formula over a finite extension, we derive a special values formula at arbitrary roots of unity that involves the discriminant of the Frobenius coinvariant N\'eron--Severi lattice.  Known bounds on the minimum discriminants of even rootless lattices then imply nontrivial lower bounds on these special values. Applying these new conditions to Kedlaya and Sutherland's list of candidate Weil polynomials of K3 surfaces over $\FF_2$ rules out over 50\% of the candidates.
\end{abstract}

\vspace*{-20pt}
\maketitle

\vspace*{-12pt}
\section*{Introduction}

Complex abelian varieties and K3 surfaces are controlled by Hodge theory, and they both enjoy well-understood period domains and Torelli theorems.  Over a finite field, our analogous understanding of the possible $\ell$-adic cohomology groups diverges between the class of abelian varieties and K3 surfaces.  Abelian varieties have a well-understood Honda--Tate theory~\cite{Honda,Tate}, in which isogeny classes are controlled by the characteristic polynomial of Frobenius, and the possible such characteristic polynomials are realized once they satisfy standard necessary conditions imposed by the Weil conjectures and crystalline cohomology.  By contrast, an analogous Honda--Tate theory for K3 surfaces over $\FF_q$, which attempts to characterize the possible Weil polynomials $P_X(t) = q\det(1-t F \;|\; H^2_{\et}(\overline X, \QQ_\ell(1)))$ of K3 surfaces $X$ over a finite field, is far from complete.  Such a theory was initiated by Elsenhans and Jahnel~\cite{elsenhans:ANTS_2010_talk,EJ} (who attribute it to Katz) and developed by Kedlaya and Sutherland~\cite{KS}, with further work in different directions by Taelman~\cite{TaL,TaO}, Ito~\cite{Ito}, and the first author with Petok~\cite{AP}.  

Besides the necessary conditions imposed by the Weil conjectures and crystalline cohomology, the possible Weil polynomials are constrained by point count nonnegativity and field extension growth~\cite{EJ,KS}, as well as by the arithmetic of special values as shown by Elsenhans and Jahnel~\cite{EJ,EJ:characteristic} using the Artin--Tate formula under base change.  An investigation of the role of noncommutative K3 surfaces~\cite{AP} suggests that these necessary conditions are not sufficient.

In this work, we construct a new class of necessary conditions on the Weil polynomial of a K3 surface, derived from techniques in the theory of sphere packing applied to the Frobenius coinvariant sublattice of the geometric N\'eron--Severi group.  Explicitly, this is equivalent to a bound on special values of the Weil polynomial at roots of unity.  For simplicity, we present the condition in the case of special values at $-1$, see Sections~\ref{sec:AT} and~\ref{sec:computation} for the general case.

\begin{maintheorem}
\label{thm:quadratic}
Let $X$ be a K3 surface over $\FF_q$ with Weil polynomial $P_X(t)$ and Picard rank $\rho$.  
Then we have the following quadratic Artin--Tate formula
\[
%
\lim_{t \to -1} \frac{P_X(t)}{(1+t)^{\rho_-} } = 
2^{\rho}
\frac{|\Br_-(X)|\, |\NS_-(X)|}{A_2(X)}
\]
 where the $(-1)$-eigenlattice of Frobenius acting on $\NS(X_{\FF_{q^2}})$ is denoted $\NS_-(X)$, a lattice of rank $\rho_-$, discrimiant $|\NS_-(X)|$, and gluing group $A_2(X)=\NS(X_{\FF_{q^2}})/\bigl(\NS(X)\oplus\NS_-(X)\bigr)$, and where $\Br_-(X) = \coker\bigl(\Br(X) \to \Br(X_{\FF_{q^2}})\bigr)$. 
Explicitly, write
\[
P_X(t) = (1-t)^{\rho} (1+t)^{\rho_-} P_2(t),
\]
with $P_2(t) \in\ZZ[t]$ such that $P_2(1) P_2(-1)\neq 0$. If $\rho_-=0$, then $P_2(-1)$ is a perfect square. If $\rho_- > 0$, then
\[ 
2^{\min(\rho,\rho^-)} \, P_2(-1) \geq \delta_{\rho_-} 
\]
where $\delta_k$ denotes the least determinant of an even positive definite lattice of rank $k$ containing no vector of norm $2$.   
\end{maintheorem}

We remark that $\rho^-$ 
coincides with the Picard rank jump from base changing to~$\FF_{q^2}$.  When the Picard rank does not jump (the case $\rho_-=0$), the statement was established by Elsenhans and Jahnel~\cite{EJ}.  The new condition is nontrivial precisely in the complementary case, when the Picard rank does jump over $\FF_{q^2}$.  Bounds on the invariant $\delta_k$, which is important in the theory of sphere packing, are summarized in Section~\ref{sec:sphere}.  In Section~\ref{sec:AT}, we present a general \emph{cyclotomic Artin--Tate formula} for the special value of the Weil polynomial at a primitive $m$th root of unity and in Section~\ref{sec:computation} combine it with the sphere packing bounds to give a new class of \emph{refined Artin--Tate} conditions generalizing Theorem~\ref{thm:quadratic}.

As an application, showing that the new conditions are quite powerful, we record their effects on the list of candidate Weil polynomials of K3 surfaces over $\FF_2$ computed by Kedlaya and Sutherland~\cite[Computation~3c]{KS}.


\begin{computation}
\label{thm:candidates}
Of the $1{,}672{,}565$ candidate Weil polynomials of K3 surfaces defined over $\FF_2$ computed by Kedlaya and Sutherland:
\begin{enumerate}\setlength{\itemsep}{5pt}
\item The quadratic refined Artin--Tate condition excludes $628{,}290$ Weil polynomials.

\item The cyclotomic refined Artin--Tate conditions (for $m > 2$) exclude $658{,}298$ Weil polynomials.

\item Together, all refined Artin--Tate conditions (for $m \geq 2$) exclude a total of $878{,}580$ Weil polynomials, approximately $52.53\%$.
\end{enumerate}
After these conditions are applied, there remain $793{,}985$ candidate Weil polynomials of K3 surfaces defined over $\FF_2$.
\end{computation}

There are two main ingredients in this work, both surrounding the discriminants of Frobenius coinvariant sublattices of the N\'eron--Severi group: 
\begin{itemize}\setlength{\itemsep}{5pt}
\item upper bounds on these discriminants are extractable from the Weil polynomial by means of special values formulas at roots of unity (the cyclotomic Artin--Tate formulas);

\item lower bounds on these discriminants follow from general sphere packing bounds on the minimum discriminants of even rootless positive definite lattices.
\end{itemize}   
Combined, these yield new nontrivial lower bounds on special values of Weil polynomials of K3 surfaces over a finite field.

To establish the cyclotomic Artin--Tate formulas, we follow an idea of Elsenhans and Jahnel~\cite{EJ} to compare the Artin--Tate formula for a surface $X$ over a finite field $\FF_q$ with its base changes to finite extensions $\FF_{q^m}$.  In \cite{EJ}, the cyclotomic Artin--Tate formulas are obtained when the Picard rank of $X$ does not jump after base change to $\FF_{q^m}$ (see Corollary~\ref{cor:EJ_higher_m}), and more work is involved in establishing the formulas and making them useful in the case when the Picard rank does jump (see Theorem~\ref{thm:cyclotomic}).  Inspired by recent work of Geisser~\cite{geisser} and the Smith theory of coinvariant lattices~\cite{boissiere_nieper-wisskirchen_sarti,boissiere_camere_mongardi}, we use calculations in Galois cohomology and the theory of gluing groups for coinvariant lattices to recast the cyclotomic Artin--Tate formula in terms of discriminants of Frobenius coinvariant lattices. Earlier celebrated work of McMullen~\cite{mcmullen:dynamics} on the dynamics of automorphisms of complex K3 surfaces, and specifically the related lattice theory developed by Gross--McMullen~\cite{gross_mcmullen} and McMullen~\cite{mcmullen:glue,mcmullen:entropy}, and further sharpened by Brandhorst~\cite{brandhorst:nonsymplectic}, relates classical resultants between cyclotomic polynomials and characteristic polynomials of orthogonal transformations to primes dividing the discriminant of the corresponding coinvariant lattice.  A related perspective is taken in work of Bayer--Fluckiger~\cite{bayer:K3,bayer:automorphisms} on characteristic polynomials of finite order orthogonal transformations.  In Section~\ref{sec:computation}, we give an exact formula for such resultants that reinterprets the  special values, which are actually norms from cyclotomic fields, making the cyclotomic Artin--Tate formulas amenable to computation.  Finally, we recall a well-known and important fact going back to Nikulin~\cite{Nikulin:finite_groups_short,Nikulin:finite_groups_long} and Kond\={o}~\cite{kondo:Niemeir}, that coinvariant N\'eron--Severi lattices of a K3 surface (after changing their sign) are positive definite with no roots.  For lattices of rank $k \leq 21$, we then use known values~\cite{ConwaySloane} (for $k \leq 8$) and bounds~\cite{CE,CdLS} (for $k > 8$) on Hermite constants to give lower bounds on the discriminants of even positive definite lattices with no roots.  This yields nontrivial lower bounds on the discriminant of Frobenius coinvariant lattices, which we can turn into lower bounds on special values of Weil polynomials.

The article is organized as follows.  Section~\ref{sec:coinvariant} contains background on the arithmetic of the coinvariant lattice associated to a finite group of isometries.  We review important consequences of the Smith theory of the coinvariant lattice~\cite{boissiere_nieper-wisskirchen_sarti,boissiere_camere_mongardi} associated to a prime order group of isometries, and prove a generalization (Lemma~\ref{lem:m}) to all cyclic groups of isometries.  Section~\ref{sec:AT} contains proofs of the quadratic and cyclotomic Artin--Tate formulas (Theorems~\ref{thm:quadratic_refined} and~\ref{thm:cyclotomic}).  Section~\ref{sec:sphere} records the fact, going back to Nikulin and Kond\={o}, that coinvariant N\'eron--Severi lattices of a K3 surface (after changing their sign) are positive definite with no roots, and then reviews existing bounds in the sphere packing literature related to the Hermite constant and minimum discriminant of an even positive definite lattice with no roots.  Section~\ref{sec:computation} ties together these ingredients.

\subsection*{Acknowledgments}

We thank Samuel Boissi\`ere, Simon Brandhorst, Phil Engel, Kiran Kedlaya, Chiara Camere, and Alessandra Sarti for helpful conversations. 
Auel was partially supported by a McLane Family Fellowship and a Walter and Constance Burke Research Initiation Award from Dartmouth College.  Both authors would like to thank ICERM for hosting the special semester on \emph{Computations on K3 Surfaces and Related Varieties} in Fall 2026, where a significant amount of this collaboration took place.

\subsection*{AI Acknowledgments}

This project emerged from a dataset (of Weil polynomials of smooth sextic and octic K3 surfaces over $\FF_2$) that was generated in the second author's undergraduate thesis~\cite{shaikh:thesis}. The second author then processed this dataset using OpenAI GPT-6 Astra to search for statistical patterns beyond those imposed by the Weil conjectures and point count nonnegativity and growth conditions (as in~\cite{KS}). Among other patterns, the model empirically discovered the quadratic refined Artin--Tate condition (and the connection with known sphere packing bounds) and suggested a proof strategy. Further prompting using OpenAI GPT-5.6, GPT-6 Astra, GPT-6 Sol, and Claude Opus 5 by the authors helped to uncover the general picture of the cyclotomic Artin--Tate formula and its related Weil polynomial conditions.  We view this work as a further indication of how large datasets of mathematical objects can lead to the discovery of surprising connections, especially as AI/ML make data analysis work much more accessible.  All text was written by the authors, who are responsible for its correctness.

\section{Invariant and coinvariant sublattices}
\label{sec:coinvariant}

The study of the coinvariant lattice of a finite group of orthogonal transformations acting on a lattice has a long history going back to work of Nikulin~\cite{Nik,Nikulin:finite_groups_short,Nikulin:finite_groups_long}, Mukai~\cite{mukai:finite_groups_k3}, and Kond\={o}~\cite{kondo:automorphisms_K3_lattices,kondo:Niemeir}, with applications to automorphisms on K3 surfaces.  Various aspects of this theory were developed in \cite{garbagnati_sarti,hashimoto}.  The Smith theory of a lattice as a module over the ring of cyclotomic integers was developed by Boissi\`ere, Nieper-Wi{\ss}kirchen, and Sarti~\cite{boissiere_nieper-wisskirchen_sarti}, which was used to study automorphisms of hyperk\"ahler manifolds, a thread that was refined in \cite{boissiere_camere_mongardi,brandhorst_cattaneo,muller}.  
We will need to generalize some aspects of the existing theory for our applications.

Let $L$ be a lattice and $G$ a finite subgroup of isometries.  We denote by $L^G$ the \linedef{invariant sublattice} and by $L_G = (L^G)^\perp$ the \linedef{coinvariant sublattice}.   Denote by $N = \sum_{g \in G} g$ the norm operator on $L$.  We recall the following, which is well-known.

\begin{lemma}
\label{lem:invar_coinvar_basics}
Let $G$ be a finite group of isometries of a lattice $L$.  We have the following:
\begin{enumerate}
\item $L^G$ is orthogonal to $L_G$ and $L^G \oplus L_G \subset L$ is of finite index,

\item $L^G \subset L$ and $L_G \subset L$ are primitive sublattices,

\item $L_G = \ker N$, 

\item $L/(L^G\oplus L_G)$ is killed by $|G|$.

\item $L/(L^G \oplus L_G) \to L^G/|G|L^G$ defined by $v \mapsto N(v) \bmod |G|L^G$ is injective.
\end{enumerate}

\end{lemma}
\begin{proof}
For (a), let $V = L \tensor \QQ$, then $\frac{1}{|G|}N : V \to V$ is the projection onto the invariant subspace $V^G = L^G \tensor \QQ$.  Since $G$ acts by isometries on $L$, we have that $\frac{1}{|G|}N$ is self-adjoint, and thus $\ker \frac{1}{|G|}N = (V^G)^\perp = V_G$.  Thus we have an orthogonal decomposition $V = V^G \oplus V_G$. For (b), the fact that $L^G = L \cap V^G$ follows since $L$ is torsion-free. 
The fact that $L_G = L \cap V_G$ follows from the fact that $V_G$ is the orthogonal complement of $V^G$. 
For (c), we then have $\ker(N : L \to L) = L \cap V_G = L_G$.  

For (d), we remark that $|G|v - N(v) = \sum_{g \in G} v-gv$ is contained in $L_G$ for every $v \in L$; indeed, it is in the kernel of $N$, using the fact that $N^2 = |G|N$.  Thus for every $v \in L$, we can write $|G|v = N(v) + (|G|v - N(v)) \in L^G \oplus L_G$ , so the finite group $L/(L^G\oplus L_G)$ is killed by $|G|$, cf.~\cite[Lemma~1.4]{mongardi_tari_wandel}.

For (e), we first note that the map is well-defined since $N$ restricts to multiplication by $|G|$ on $L^G$ and to zero on $L_G$.  Now assume that $N(v) = |G| w$ for some $w \in L^G$. Then $N(v-w) = 0$, which implies that $v-w \in L_G$, so that $v = w + (v-w) \in L^G \oplus L_G$.  Thus the map is injective. 
\end{proof}

The group $L/(L^G \oplus L_G)$ is often referred to as a \emph{gluing group}, cf.\ \cite[Chapter~4.3]{ConwaySloane}.

\begin{remark}
\label{rem:eigenlattice}
When $G$ is a cyclic group of order $m$ generated by $\sigma$, then $L \tensor \CC$ decomposes into eigenspaces corresponding to eigenvalues that are $m$th roots of unity.  Then $L^G \tensor \CC$ is the $(+1)$-eigenspace and $L_G \tensor \CC$ is the direct sum of all other eigenspaces. 

In particular, when $G$ is a cyclic group of order $2$, generated by an involution $\sigma$, then $L^G = L_+$ is the $(+1)$-eigenlattice and $L_G = L_-$ is the $(-1)$-eigenlattice of $\sigma$.  
\end{remark}

When $G$ is a cyclic group of order $m$ generated by $\sigma$, we have that $N(1-\sigma) = 1 - \sigma^m = 0$, hence that $1-\sigma : L \to L$ has kernel $L^G$ and image contained in $L_G = \ker N$.  Via restriction, we can then consider the endomorphism $1-\sigma : L_G \to L_G$.

\begin{lemma}
\label{lem:det}
Let $G = \langle \sigma \rangle$ be a cyclic group of order $m > 1$ of isometries of a lattice $L$, and denote by $A = L/(L^G \oplus L_G)$.  Then $1-\sigma : L_G \to L_G$ has full rank and we have
\[
[L_G : (1-\sigma) L_G] = \det\bigl( 1-\sigma \, | \, L_G \otimes \QQ \bigr)
\]
and
\[
|H^1(G,L)| = [L_G : (1-\sigma) L] = \frac{1}{|A|}\det\bigl( 1-\sigma \, | \, L_G \otimes \QQ \bigr).
\]
\end{lemma}
\begin{proof}
The fact that $1-\sigma : L_G \to L_G$ has full rank (for $m > 1$) follows from Lemma~\ref{lem:invar_coinvar_basics}.  In particular, $1-\sigma$ is injective when restricted to $L_G$.  Since $L_G$ is a free abelian group, we get the first equality, since we can compute the determinant of $1-\sigma$ on a $\ZZ$-basis of $L_G$ or on a $\QQ$-basis of $L_G \otimes \QQ$.

For the second equality, we have $H^1(G,L) \isom L_G/(1-\sigma)L$ by the cohomology of cyclic groups.  To compare the index of $(1-\sigma)L$ and of $(1-\sigma)L_G$ in $L_G$, we consider the commutative diagram
\[
\xymatrix{
0 \ar[r] & L_G \ar[r]^{1-\sigma} \ar[d] & L_G \ar[r] \ar@{=}[d] & L_G/(1-\sigma)L_G \ar[r] \ar@{->>}[d] & 0 \\
0 \ar[r] & L/L^G \ar[r]^{1-\sigma} & L_G \ar[r] & L_G/(1-\sigma)L \ar[r] & 0
}
\]
which implies by the snake lemma that the kernel of the right vertical map is isomorphic to $(L/L^G)/L_G \isom L/(L^G\oplus L_G)$ since $L^G \cap L_G = 0$.  Thus $[L_G : (1-\sigma)L_G] = |A| [L_G : (1-\sigma)L]$.  This proves the second equality.
\end{proof}

\begin{remark}
\label{rem:2^r}
When $\sigma$ is an involution on a lattice $L$ with eigensublattices $L^+$ and $L^-$ corresponding to eigenvalue $+1$ and $-1$, respectively, then by Lemma~\ref{lem:det} (and an immediate generalization), we have that $[L : L_+ \oplus L_-]\, [L_\pm : (1\pm\sigma)L] = 2^{r_\pm}$, and in particular, that the $\FF_2$-dimension of $L/(L_+\oplus L_-)$ is bounded above by $\min(r_+,r_-)$.  This fact is contained in Nikulin's work on involutions on lattices and 2-elementary lattices. 
\end{remark}

More generally, when $G$ has prime order, one has the following restrictions on the coinvariant lattice and the ``gluing group'' $L/(L^G \oplus L_G)$, cf.\ \cite[Lemma~1.8]{mongardi_tari_wandel}.

\begin{lemma}
\label{lem:mprime}
Let $G$ be a group of prime order $p$ acting on a lattice $L$.  Then there exists $a_p \geq 0$ such that
\[
\rk L_G = a_p (p-1).
\]
Writing $[L : L^G\oplus L_G] =p^a$, we have 
\[
a \leq \min\{v_p(|L^G|),v_p(|L_G|),a_p\}
\]
where $|L^G|$ and $|L_G|$ denote the absolute value of the discriminants of the invariant and coinvariant sublattices.
\end{lemma}

We provide a generalization of Lemma~\ref{lem:mprime} to arbitrary cyclic group actions.

\begin{lemma}
\label{lem:m}
Let $G$ be a cyclic group of order $m$ acting on a lattice $L$.  Then for each $d \mid m$, there exists $a_d \geq 0$ such that
\[
\rk L_G = \sum_{d \mid m,\, d>1} a_d \varphi(d).
\]
Writing $A = L/(L^G\oplus L_G)$, then for every prime $p \mid m$, we have 
\[
v_p(|A|) \leq \min\bigl\{v_p(|L^G|),v_p(|L_G|),b_p(L_G)\}
\]
where $|L^G|$ and $|L_G|$ denote the absolute value of the discriminants of the invariant and coinvariant sublattices, and where we define
\[
b_p(L_G) = \sum_{i=1}^{v_p(m)} a_{p^i}
\]
Also, we have that
\[
\ell(A_p) \leq (\rk L_G)/(p-1)
\]
where $\ell(A_p)$ denotes the number of cyclic factors of the $p$-primary part of $A$. 
\end{lemma}
\begin{proof}
The action of $G = \langle \sigma \rangle$ extends to a representation of $G$ on $V_G=L_G\tensor \QQ$, yielding a decomposition
\[
V_G \isom \bigoplus_{d \mid m,\, d>1} \rho_d^{\oplus a_d}
\]
into nontrivial (since $L_G$ is orthogonal to $L^G$) irreducible representations $\rho_d$ of $G$, for some multiplicities $a_d \geq 0$.  Here, $\rho_d = \QQ(\zeta_d)$ with $\sigma$ acting by multiplication by $\zeta_d$.  By Lemma~\ref{lem:det}, we have that
\[
|H^1(G,L)| = [L_G:(1-\sigma)L] = \frac{1}{|A|}\det\bigl( 1-\sigma \, | \, V_G \bigr)
\]
which, since $|H^1(G,L)|$ is an integer, shows that $|A|$ divides $\det(1-\sigma \, | \, V_G)$.  
Since $\sigma$ acts by multiplication by $\zeta_d$ on each summand $\rho_d$ of $V_G$, the characteristic polynomial of $\sigma$ on $V_G$ is $\prod_{d \mid m, d > 1} \Phi_d(T)^{a_d}$.  Thus
\[
\det\bigl( 1-\sigma \, | \, V_G \bigr) = \prod_{d \mid m,\, d > 1} \Phi_d(1)^{a_d}.
\]
However, by a formula first discovered by E.~Lehmer~\cite[Theorem~3]{lehmer:cyclotomic} and later in \cite{apostol,Diederichsen}, cf.\ \cite[Theorem~3.14]{brandhorst:nonsymplectic}, \cite[Lemma~4.1]{bayer:K3}, we have
\begin{equation}
\label{eq:cyclotomic1}
\Phi_d(1) = 
\begin{cases}
0 & \text{if $d=1$} \\
p & \text{if $d = p^k$ is a prime power with $k > 0$}\\
1 & \text{if $d$ is not a prime power}
\end{cases}
\end{equation}
which for $d > 1$ coincides with $e^{\Lambda(d)}$ where $\Lambda$ is the von Mangoldt function from analytic number theory. 
Thus we have that
\begin{equation}
\label{eq:det}
\det(1-\sigma \, | \, V_G) = \prod_{p \mid m} p^{b_p(L_G)}
\end{equation}
and so we arrive at the desired inequality of $p$-adic valuations
\[
v_p(|A|) \leq b_p(L_G) 
\]
Since $L^G$ and $L_G$ are primitive in $L$, we have that $A = L/(L^G\oplus L_G)$ embeds into the direct sum $A_{L^G} \oplus A_{L_G}$ of discriminant groups, and projecting to each factor yields embeddings $A \subset A_{L^G}$ and $A \subset A_{L_G}$, see \cite[Section~1.5]{Nik}.  This gives
\[
v_p(|A|) \leq \min\{v_p(|L^G|),v_p(|L_G|)\},
\]
which we combine with the previous inequality.  Finally, to bound the number of cyclic factors of the $p$-primary part of $A$, we have that
\[
\ell(A_p) \leq b_p(L_G) 
\leq \frac{1}{p-1} \sum_{i=1}^{v_p(m)} a_{p^i}\varphi(p^i) \leq \frac{\rk L_G}{p-1}
\]
as desired.
\end{proof}

\section{Artin--Tate special value formulas at roots of unity}
\label{sec:AT}

Let $X$ be a K3 surface over a finite field $\FF_q$.  We write $\rho = \rk \NS(X)$ for the Picard rank of $X$ and $|\NS(X)|$ for the absolute value of the discriminant of the N\'eron--Severi lattice of $X$.  We write $P_X(t) = q\det(1-t F \;|\; H^2_{\et}(\overline X, \QQ_\ell(1)))$ for the normalized $L$-polynomial of Frobenius acting on $\ell$-adic cohomology, which we call the \linedef{Weil polynomial} of $X$.  Then $P_X(t) \in \ZZ[t]$ has degree 22, all roots of absolute value $1$, and $P_X(0)=q$.

For a rational function $P(t) \in \Qbar(t)$ and $a \in \Qbar$ nonzero, we define the \linedef{special value} of $P$ at $a$ to be
\begin{equation}
\label{eq:special_value}
P^*(a) = 
\lim_{t \to a} \frac{P(t)}{(1-t/a)^r}
=
\lim_{t \to 1} \frac{P(at)}{(1-t)^r}
\end{equation}
where $r$ is the order of vanishing of $P(t)$ at $t=a$.

We recall the Artin--Tate formula, which is a special value formula for the Weil polynomial at $t=1$,
\begin{equation}
\label{eq:AT}    
P_X^*(1) = \lim_{t\to 1}
\frac{P_{X}(t)}{(1-t)^{\rho}}
=
|\Br(X)|\, |\NS(X)|
\end{equation}
whose validity has been established by the proof of the Tate conjecture for K3 surfaces over a finite field~\cite{NygaardOgus,Charles,MadapusiPera,KimMadapusiPera,MPerr,IIK}.  We will be interested in comparing the Artin--Tate formula for various base changes of $X$.

To that end, for any $m \geq 1$ and any smooth projective surface $X$ over $\FF_q$, we write $X_m = X_{\FF_{q^m}}$ for the base change of $X$ to $\FF_{q^m}$ and $\rho_m = \rk \NS(X_m)$.  We denote by $\NS_m(X) = \NS(X)^\perp \cap \NS(X_m)$ the coinvariant lattice of the (cyclic subgroup generated by the) Frobenius of $X$ acting on $\NS(X_m)$, and by $A_m(X) = \NS(X_m)/(\NS(X)\oplus \NS_m(X))$ the corresponding gluing group.  We also denote by $\Br_m(X) = \coker\bigl(\Br(X) \to \Br(X_m)\bigr)$.  We take $\NS_1(X)=0$ (which has discriminant 1), $A_1(X)=0$, and $\Br_1(X)=0$.

Since the action of Frobenius on $H^2_{\et}(\overline{X}_m, \QQ_\ell(1))$ is the $m$th power of the action on $H^2_{\et}(\overline X, \QQ_\ell(1))$, the formula for the characteristic polynomial of the $m$th power of an operator immediately yields the following factorization formula for the Weil polynomial
\begin{equation}
\label{eq:AT_base_change}
P_{X_m}(t^m)
=
\prod_{\zeta^m=1}P_X(\zeta t).
\end{equation}
Here, we remind the reader that $P_{X_m} = q^m \det(1-tF^m \, | \, H^2_{\et}(\overline X, \QQ_\ell(1)))$.  We remark that, via the Weil conjectures, this factorization can also be derived from the generating function description of the zeta function, cf.\ \cite[Proposition~2.14]{mustata:zeta_book}.

To state our special value formulas, we will need the following standard application of the Leray spectral sequence to the restriction map on the Brauer group, cf.\ \cite[Section~2]{geisser}.  For the convenience of the reader, we provide a proof.

\begin{lemma}
\label{lem:descent}
Let $X$ be a smooth projective variety over a finite field $\FF_q$.  Then for every $m \geq 1$ we have
\[
 \ker\bigl(\Br(X) \to \Br(X_m)\bigr) \isom H^1(G,\NS(X_m)).
\]
where $G=\Gal(\FF_{q^m}/\FF_q)$.
\end{lemma}
\begin{proof}
More generally, for any geometrically integral scheme $X$ proper over a field $k$ and a finite $G$-Galois extension field $l/k$, the Leray spectral sequence associated to the structure morphism $X_l \to X$ and the multiplicative group $\Gm$ (i.e., the Hochschild--Serre spectral sequence)
\[
H^i(G,H^j(X_l,\Gm)) \Longrightarrow H^{i+j}(X,\Gm)
\]
has sequence of low degree terms
\begin{align*}
0 \longrightarrow {} & H^1(G,l^\times) \longrightarrow  \Pic(X) \longrightarrow \Pic(X_l)^G \longrightarrow H^2(G,l^\times) \\ 
 \longrightarrow {} & \ker\bigl(\Br(X)\to \Br(X_l)\bigr) \longrightarrow H^1(G,\Pic(X_l))
 \longrightarrow H^3(G,l^\times).
\end{align*}
By Hilbert's Theorem 90, we have that $H^1(G,l^\times)=0$.  If we now assume that $l/k$ is a finite cyclic extension, then periodicity for cohomology of a cyclic group gives $H^3(G,l^\times) = H^1(G,l^\times) = 0$ and $H^2(G,l^\times) = k^\times/N_{l/k}l^\times$.  Finally, when $k=\FF_q$ is a finite field and $l = \FF_{q^m}$, this latter group is trivial.  In this case, we deduce that the map
\[
\ker\bigl(\Br(X)\to \Br(X_m)\bigr) \longrightarrow H^1(G,\Pic(X_m))
\]
is an isomorphism.  Now consider the short exact sequence
\[
0 \to \Pic^0(X_m) \to \Pic(X_m) \to \NS(X_m) \to 0
\]
of $G$-modules. Lang's theorem on vanishing of higher Galois cohomology of abelian varieties over a finite field gives $H^i(G,\Pic^0(X_m))=0$ for all $i \geq 1$, and then the long exact sequence yields an isomorphism $H^i(G,\Pic(X_m)) \to H^i(G,\NS(X_m))$ for all $i \geq 1$.
\end{proof}

\medskip

As a warm-up for our special value calculations, we first give a formula, which is implicit in \cite[Theorem~12]{EJ}, for the special value of the Weil polynomial of a K3 surface at $-1$.  To that end, we remark that when $m=2$, we have that $\NS_2(X) = \NS_-(X)$
is the $(-1)$-eigenlattice of Frobenius acting on $\NS(X_2)$, as in Remark~\ref{rem:eigenlattice}.  In analogy, we also write $\Br_-(X) = \Br_2(X)$.

\begin{theorem}
\label{thm:quadratic_refined}
Let $X$ be a K3 surface over $\FF_q$ with Weil polynomial $P_X(t)$.  Then we have the following quadratic Artin--Tate formula
\[
P_X^*(-1) =
2^{\rho_2}
\frac{|\Br(X_2)|}{|\Br(X)|}
\frac{|\NS(X_2)|}{|\NS(X)|}
= 
2^{\rho} \, 
\frac{|\Br_-(X)| \, |\NS_-(X)|}{|A_2(X)|}.
\] 
\end{theorem}
\begin{proof}
By \eqref{eq:AT_base_change}, we have $P_{X_2}(t^2) = P_X(t)\, P_X(-t)$.  Dividing both sides by $(1-t^2)^{\rho_2}$, we have
\[
\frac{P_{X_2}(t^2)}{(1-t^2)^{\rho_2}} = (1+t)^{-\rho_2}\frac{P_X(t)}{(1-t)^\rho} \frac{P_X(-t)}{(1-t)^{\rho^-}}
\]
where we write $\rho_2 = \rho + \rho_-$, where $\rho_- = \rk \NS_-(X)$. Taking limits as $t \to 1$, we have
\[
|\Br(X_2)| \, |\NS(X_2)| = 2^{-\rho_2} \, |\Br(X)| \, |\NS(X)| \, P_X^*(-1)
\]
where the left hand side follows from the Artin--Tate formula~\eqref{eq:AT} for $X_2$, and the right hand side follows from the Artin--Tate formula for $X$ and  (by a change of variables $t \mapsto -t$) the special value at $t=-1$, which yields the first desired formula for $P_X^*(-1)$.

For the second formula, the exact sequence
\[
0 \to \ker\bigl( \Br(X) \to \Br(X_2) \bigr) \to \Br(X) \to \Br(X_2) \to \Br_-(X) \to 0
\]
combined with Lemmas~\ref{lem:det} and~\ref{lem:descent}, yields
\[
\frac{|\Br(X_2)|}{|\Br(X)|} = 
\frac{|\Br_-(X)|}{|H^1(G, \NS(X_2))|} = \frac{|\Br_-(X)||A_2(X)|}{2^{\rho^-}}
\]
Finally, by comparing discriminants, we have 
\[
|\NS(X)|\, |\NS_-(X)| = |\NS(X) \oplus \NS_-(X)| = |A_2(X)|^2 \, |\NS(X_2)|
\]
from which the stated result follows. 
\end{proof}

As a corollary, we arrive at (a refinement of) the condition of Elsenhans and Jahnel~\cite{EJ}.

\begin{corollary}
\label{cor:EJ_formula}
Let $X$ be a K3 surface over $\FF_q$ with Weil polynomial $P_X(t)$. Assume that $\rho = \rho_2$.  Then $\NS(X) \to \NS(X_2)$ is an isomorphism and $\Br(X) \to \Br(X_2)$ is injective, hence the quadratic Artin--Tate formula reduces to
\[
P_X^*(-1) = 2^{\rho} \, |\Br_-(X)|.
\]
where $\Br_-(X)$ has (integral) square order.
\end{corollary}
\begin{proof}
The base change map $\NS(X) \to \NS(X_2)$ is injective by a restriction-corestriction argument.  Hence the $(-1)$-eigenspace of Frobenius acting on $\NS(X_2)\tensor \QQ$ is trivial.  And since $\NS(X_2)$ is torsion-free, the $(-1)$-eigenlattice $\NS_-(X)$ must be trivial as well.  Now, for $v \in \NS(X_2)$ we have that $F(v)-v$ is contained in $\NS_-(X)$, hence must be zero, showing that $v \in \NS(X_2)^F = \NS(X)$, where this last equality follows from the Leray spectral sequence associated to $X_2 \to X$ and the sheaf $\Gm$, the triviality of the $\Pic^0$ scheme of a K3 surface, and the fact that the Brauer group of a finite field is trivial.  Thus the base change map $\NS(X) \to \NS(X_2)$ is surjective as well.

By Lemma~\ref{lem:descent}, the kernel of the base change map $\Br(X) \to \Br(X_2)$ is isomorphic to $H^1(G,\NS(X_2))$, where $G$ is generated by Frobenius.  This latter group is trivial as $\NS(X_2)$ is a free abelian group with trivial $G$-action.  Hence the base change map on Brauer groups is injective, and the cokernel $\Br_-(X)$ has order $|\Br(X_2)|/|\Br(X)|$, which is an integer that is a quotient of squares (as the Brauer group of a surface over a finite field has square order), hence is itself a perfect square.
\end{proof}

Now, we present the following generalization of Theorem~\ref{thm:quadratic_refined} to give an Artin--Tate formula for special values at higher degree roots of unity.  To that end, for a variety $X$ over $\FF_q$ define $b_1(X)=1$ and for $m \geq 2$ define
\[
b_m(X) = \prod_{p \mid m} p^{b_p(\NS_m(X))}
\]
where the product is over primes $p$ dividing $m$ and where the invariant $b_p$ of a coinvariant lattice is introduced in the statement of Lemma~\ref{lem:m}.  For example, when $m=p$ is prime  Lemma~\ref{lem:mprime} gives $b_p(X) = p^{\rk \NS_p(X)/(p-1)}$.  In particular, in the context of Theorem~\ref{thm:quadratic_refined}, $b_2(X) = 2^{\rho_-}$.

\begin{theorem}
\label{thm:cyclotomic}
Let $X$ be a K3 surface over $\FF_q$ with Weil polynomial $P_X(t)$.  Then for every $m \geq 2$ we have the following cyclotomic Artin--Tate formula
\begin{align*}
N_{\QQ(\zeta_m)/\QQ}(P_X^*(\zeta_m))
&{} = 
\prod_{d\mid m}
\bigl(
d^{\rho_d}\, |\Br(X_d)|\, |\NS(X_d)|
\bigr)^{\mu(m/d)} \\
&{} =
\prod_{d\mid m}
\left(
\frac{d^{\rho_d}\, |\Br_d(X)|\, |\NS_d(X)|}
{|A_d(X)|\, b_d(X)}
\right)^{\mu(m/d)}.
\end{align*}
\end{theorem}
\begin{proof}
We proceed as in the proof of Theorem~\ref{thm:quadratic_refined}, with additional elaboration required in the cyclotomic case.  By \eqref{eq:AT_base_change}, we have $P_{X_m}(t^m) = \prod_{\zeta^m=1} P_X(\zeta t)$.  Writing $\rho_\zeta$ for the algebraic multiplicity of $\zeta$ as an eigenvalue of Frobenius acting on $\NS(X_m)$, we have that $\rho_m = \sum_{\zeta^m=1} \rho_\zeta$.  Dividing both sides by $(1-t^m)^{\rho_m}$, we have
\[
\frac{P_{X_m}(t^m)}{(1-t^m)^{\rho_m}} = \frac{1}{(1+t+\dotsm+t^{m-1})^{\rho_m}}\prod_{\zeta^m=1}\frac{P_X(\zeta t)}{(1-t)^{\rho_\zeta}}
\]
and taking limits as $t \to 1$ and using \eqref{eq:special_value}, we have
\[
m^{\rho_m}\, P_{X_m}^*(1) = \prod_{\zeta^m=1} P_X^*(\zeta) = \prod_{d \mid m} N_{\QQ(\zeta_d)/\QQ}(P_X^*(\zeta_d))
\]
where we use the fact that $P_X(t)$ has rational coefficients (in fact, integral coefficients, but we don't need that here). M\"obius inversion and the Artin--Tate formula~\eqref{eq:AT} for $X_m$ then give
\[
N_{\QQ(\zeta_m)/\QQ}(P_X^*(\zeta_m)) = \prod_{d \mid m} \bigl(d^{\rho_d}\, |\Br(X_d)|\, |\NS(X_d)|\bigr)^{\mu(m/d)}
\]
which is the first claimed equality in the formula.

For the second equality, the exact sequence
\[
0 \to \ker\bigl( \Br(X) \to \Br(X_d) \bigr) \to \Br(X) \to \Br(X_d) \to \Br_d(X) \to 0
\]
combined with Lemmas~\ref{lem:det} and~\ref{lem:descent} yields
\[
|\Br(X_d)| = 
|\Br(X)|\, 
\frac{|\Br_d(X)|}{|H^1(G, \NS(X_d))|}.
\]
By comparing discriminants, we have 
\[
|\NS(X)|\, |\NS_d(X)| = |\NS(X) \oplus \NS_d(X)| = |A_d(X)|^2 \, |\NS(X_d)|
\]
Thus, using Lemmas~\ref{lem:det} and ~\ref{lem:descent}, we have 
\[
|\Br(X_d)|\, |\NS(X_d)| 
= 
|\Br(X)|\, |\NS(X)|\, 
\frac{|\Br_d(X)|\, |\NS_d(X)|}{b_d\, |A_d(X)|}
\]
where we temporarily set $b_d = |\det(1-F \, | \, \NS_d(X)\tensor \QQ)|$.  By \eqref{eq:det}, we compute that
\[
b_d = \det(1-F \, | \, \NS_d(X)\tensor \QQ)
= \prod_{p|m} p^{b_p(\NS_d(X))} = b_d(X).
\]
Finally, since $\prod_{d \mid m} c^{\mu(m/d)} = 1$ for any $c \neq 0$ and any $m \geq 2$, we have that
\begin{align*}
N_{\QQ(\zeta_m)/\QQ}(P_X^*(\zeta_m)) &{} = \prod_{d \mid m} \left(\frac{d^{\rho_d}\, |\Br_d(X)|\, |\NS_d(X)|}{b_d(X)\, |A_d(X)|} \right)^{\mu(m/d)}
\end{align*}
which is the second stated equality.  
\end{proof}

\begin{remark}
\label{rem:P_1}
Considering a different normalization of the special value leads to some simplifications. To this end, write 
\[
P_X(t) = (1-t)^\rho \, P_1(t).
\]
Then $P_X^*(\zeta_m) = (1-\zeta_m)^\rho P_1^*(\zeta_m)$
 implies that
\[
N_{\QQ(\zeta_m)/\QQ}(P_X^*(\zeta_m)) = \Phi_m(1)^\rho N_{\QQ(\zeta_d)/\QQ}(P_1^*(\zeta_d))
\]
By \eqref{eq:cyclotomic1}, we have $\prod_{d \mid m,\, d > 1} \Phi_d(1) = m$, which implies that $\prod_{d \mid m} d^{\mu(m/d)} = \Phi_m(1)$ for $m \geq 2$ by M\"obius inversion.  Hence Theorem~\ref{thm:cyclotomic} can be restated as
\begin{align}
\label{eq:cyclotomic_P_1}
\begin{split}
N_{\QQ(\zeta_m)/\QQ}(P_1^*(\zeta_m))
&{} = 
\prod_{d\mid m}
\bigl(
d^{\rho_d-\rho}\, |\Br(X_d)|\, |\NS(X_d)|
\bigr)^{\mu(m/d)} \\
&{} =
\prod_{d\mid m}
\left(
\frac{d^{\rho_d-\rho}\, |\Br_d(X)|\, |\NS_d(X)|}
{|A_d(X)|\, b_d(X)}
\right)^{\mu(m/d)}
\end{split}
\end{align}
In the setting of Corollary~\ref{cor:EJ_formula}, when $\rho = \rho_2$ we know that $(1+t)$ does not divide $P_X(t)$, and we have the strikingly simple formula
\[
P_1(-1) = P_1^*(-1) = |\Br^-(X)|
\]
for the order of the growth of the Brauer group under base change to $\FF_{q^2}$.
\end{remark}

We have the following generalization of the construction of Elsenhans and Jahnel~\cite{EJ} to the higher cyclotomic case, which is analogous to (and proved in the same way as) Corollary~\ref{cor:EJ_formula}.

\begin{corollary}
\label{cor:EJ_higher_m}
Let $X$ be a K3 surface over $\FF_q$ with Weil polynomial $P_X(t)$. For some $m \geq 2$, assume that $\rho = \rho_m$.  Then $\NS(X) \to \NS(X_m)$ is an isomorphism and $\Br(X) \to \Br(X_m)$ is injective, and the cyclotomic Artin--Tate formula reduces to
\[
N_{\QQ(\zeta_m)/\QQ}(P_X^*(\zeta_m)) = \Phi_m(1)^\rho \, \prod_{d \mid m} |\Br_d(X)|^{\mu(m/d)}
\]
or equivalently
\[
N_{\QQ(\zeta_m)/\QQ}(P_1^*(\zeta_m)) = \prod_{d \mid m} |\Br_d(X)|^{\mu(m/d)}
\]
where the alternating product of Brauer cokernels is an integral square. 
\end{corollary}
\begin{proof}
The statements about the base change maps on the N\'eron--Severi and Brauer groups follow exactly as in the proof of Corollary~\ref{cor:EJ_formula}.  In this case, $\rho = \rho_d$ and $|\NS_d(X)|=|A_d(X)|=b_d = 1$ for all $d \mid m$, hence we arrive at the desired formulas by applying \eqref{eq:cyclotomic1} and Remark~\ref{rem:P_1}.  

Finally, we consider the alternating product of Brauer cokernels.  Since the order of the Brauer group of a K3 surface over a finite field is a square and $\Br(X) \to \Br(X_d)$ is injective for all $d \mid m$, we have that $\Br_d(X)$ has square order.  Thus the alternating product is a rational square.  Since $P_1(t)$ has integer coefficients, $N_{\QQ(\zeta_m)/\QQ}(P_1^*(\zeta_m))$ is an algebraic integer, and thus so is the alternating product of Brauer cokernels.  A rational square that is an algebraic integer is an integer square, yielding the final claim.
\end{proof}

\begin{remark}
\label{rem:EJ_higher}
We now make a comparison with the construction of \cite[Section~5]{EJ}.  Elsenhans and Jahnel define an invariant $E^{(m)}$ associated to any degree 22 Weil polynomial with roots of absolute value $q$.  Translating to our Weil polynomial convention (roots of absolute value 1), we have that
\[
E^{(m)} = q^{m-1} \, \prod_{\alpha \neq 1} \frac{1-\alpha^m}{1-\alpha}
\]
where the product is over all inverse roots $\alpha$ of a Weil polynomial $P$.  When $P=P_X$ for a K3 surface $X$ over $\FF_q$, we see that $E^{(m)}=0$ if and only if $\rho_m > \rho$, cf.\ \cite[Remark~13]{EJ}.  In the case $\rho=\rho_m$, we have that
\[
E^{(m)} = \prod_{\zeta^m=1,\, \zeta \neq 1} \!P_1(\zeta) 
= 
\prod_{d \mid m, \, d>1} \!N_{\QQ(\zeta_d)/\QQ}(P_1(\zeta)) 
=
\frac{1}{m^\rho} \prod_{d \mid m, \, d>1} \!N_{\QQ(\zeta_d)/\QQ}(P_X^*(\zeta_d)). 
\]
The first equation follows from the fact that $(1-\alpha^m)/(1-\alpha) = \prod_{\zeta^m=1,\, \zeta\neq 1} (1-\zeta\alpha)$.  The final equation follows from Remark~\ref{rem:P_1}. On the other hand, still assuming $\rho=\rho_m$, we have
\[
E^{(m)} = \frac{P_{X_m}^*(1)}{P_X^*(1)} = \frac{|\Br(X_m)| \, |\NS(X_m)|}{|\Br(X)| \, |\NS(X)|} = |\Br_m(X)|
\]
using the classical Artin--Tate formula and the fact that $|\Br_m(X)| = |\Br(X_m)|/|\Br(X)|$.  This gives a refinement of the fact that $E^{(m)}$ is a square, see~\cite[Theorem~12]{EJ}.
\end{remark}

\section{Sphere packing bounds}
\label{sec:sphere}

In this section, we survey the existing literature for lower bounds on the discriminants of lattices arising in the theory of sphere packing.  Such bounds are useful for K3 surfaces over a finite field because our cyclotomic Artin--Tate formulas give upper bounds for the discriminants of Frobenius coinvariant lattices in terms of special values of the Weil polynomial. The following result, which in spirit goes back to Nikulin~\cite{Nikulin:finite_groups_short,Nikulin:finite_groups_long}, cf.\ Kond\={o}~\cite[Proposition~2]{kondo:Niemeir} motivates the entire approach. 

\begin{lemma}
\label{lem:rootless}
Let $X$ be a smooth K3 surface over a field $k$, and let $l/k$ be an extension.  Then the sublattice $N = \NS(X)^\perp \cap \NS(X_l)$ is even, negative definite, and has no $(-2)$-classes.
\end{lemma}

\begin{proof}
Let $H$ be an ample divisor on $X$ defined over $k$.  It is a standard fact that $\NS(X_l)$ is even, hence so is any sublattice. Since $N \subset H^\perp$, the Hodge index theorem implies that $N$ is negative definite. We now prove that this sublattice is rootless. In the case where $N=0$ the statement is trivial. Assume that $C \in N$ satisfies $C^2=-2$, then either $C$ or $-C$ is effective (cf.\ \cite[Section~2.1.5]{huybrechts:K3_book}), in which case $H.C=0$ is not possible by Nakai--Moishezon--Kleiman (cf.\ \cite[Theorem~8.1.2]{huybrechts:K3_book}). Hence $N$ has no $(-2)$-classes.
\end{proof}

Note that the conclusion of Lemma~\ref{lem:rootless} does not imply that $X$ cannot acquire $(-2)$-classes over an extension, but just that any $(-2)$-class cannot be orthogonal to $\NS(X)$.

\begin{remark}
For $l/k$ a finite $G$-Galois extension, note that, in the notation of Lemma~\ref{lem:rootless}, $N = N_G$ is the coinvariant lattice of the induced action of $G$ on $\NS(X_l)$.  Hence, in spirit, Lemma~\ref{lem:rootless} is an analogue of Nikulin~\cite[Lemma~4.2]{Nikulin:finite_groups_long}.
\end{remark}

Now, let $L = N(-1)$ be obtained by changing the sign of the intersection form. Then $L$ is even, positive definite, and \emph{rootless}, i.e., contains no vectors of norm $2$. In particular, the norm of any vector in $L$ is at least $4$. Positive definite lattices satisfying this property (even or not) have been extensively studied as they describe integral lattice sphere packings.

More specifically, for a positive definite lattice $L$ of rank $k$, define its \emph{minimum} $\mu(L)$ to be the least norm of a non-zero vector in $L$ and its \emph{discriminant} $D(L)$ to be the determinant of a Gram matrix of an integral basis of $L$. Then the minimum distance between distinct lattice points is $\sqrt{\mu(L)}$ and the volume of a fundamental domain for $L$ is given by $\sqrt{D(L)}$. Consequently, balls of radius $\frac{1}{2}\sqrt{\mu(L)}$ centered at the lattice points form a packing of density 
\[ 
\Delta(L) = \frac{v_k}{2^k} \frac{\mu(L)^{k/2}}{D(L)^{1/2}} 
\]
where $v_k$ is the volume of the unit ball in $\mathbb R^k$. 

The classical problem of finding the densest lattice sphere packing in dimension $k$ is equivalent to determining the \emph{Hermite constant}
\[ 
\gamma_k=\sup_{L} \frac{\mu(L)}{D(L)^{1/k}}, 
\]
where $L$ ranges over all positive definite lattices of rank $k$. Indeed, let $\Delta_k$ be the largest density of a lattice sphere packing in dimension $k$. Then we have that 
\[ 
\Delta_k = \frac{v_k}{2^k}\gamma_k^{k/2}. 
\]
For our purposes, we can use Hermite's inequality $\mu(L) \leq \gamma_k D(L)^{1/k}$ to provide a lower bound on $D(L)$ in terms of the minimum $\mu(L)$. Let $\delta_k$ denote the least discriminant of a positive definite even rootless lattice of rank $k$, hence with minimum at least 4. Then we have
\[ 
D(L) \geq \delta_k \geq \left(\frac{4}{\gamma_k}\right)^{\!\!k}.
\]
The exact Hermite constants are known only in rank $k \leq 8$, see~\cite[Chapter~1]{ConwaySloane}, and $k=24$~\cite{CohnKumar}. In rank up to $8$, we have 
\[
\renewcommand{\arraystretch}{1.5}
\begin{array}{c|rrrrrrrr}
k & 1 & 2 & 3 & 4 & 5 & 6 & 7 & 8 \\\hline
\gamma_k^k & 1 & \frac{4}{3} & 2 & 4 & 8 & \frac{64}{3} & 64 & 256 \\\hline
\delta_k & 4 & 12 & 32 & 64 & 128 & 192 & 256 & 256
\end{array}
\]
which are attained by the scaled even root lattices $A_1(2)$, $A_2(2)$, $A_3(2)$, $D_4(2)$, $D_5(2)$, $E_6(2)$, $E_7(2)$, and $E_8(2)$~\cite[Chapter~1]{ConwaySloane}.  In rank $24$, the Hermite constant is attained by the Leech lattice.  In the remaining ranks up to $21$, we can appeal to numerical upper bounds on the Hermite constant established by Cohn and Elkies~\cite{CE}, and updated by Cohn, de Laat, and Salmon~\cite{CdLS} to provide lower bounds on lattice discriminants. 
\[
\renewcommand{\arraystretch}{1.5}
\begin{array}{c|rrrrrrrrrrrrr}
k & 9 & 10 & 11 & 12 & 13 & 14 & 15 & 16 & 17 & 18 & 19 & 20 & 21 \\\hline
\delta_k\ge & 298 & 317 & 312 & 281 & 230 & 177 & 128 & 89 & 58 & 37 & 24 & 13 & 8
\end{array}
\]

\begin{remark}
By Lemma~\ref{lem:rootless}, the negatives of coinvariant N\'eron--Severi lattices of K3 surfaces are positive definite even rootless lattices, which is why we restrict our attention to even lattices.  Since the Hermite constant is invariant under scaling, we can always assume that lattices are even for the purposes of bounding the Hermite constant.  However, the minimum discriminant of positive definite even rootless lattices can, a priori, be larger than the minimum discriminant of all positive definite lattices with minimum at least~$4$.  Thus the invariant $\delta_k$ that we care about may satisfy even stronger bounds than those achievable via an analysis of the Hermite constant, see \cite{pointless} for example.  When $k \leq 8$ and $k=24$, in fact the even lattices achieving the bound (listed above) have minimum~$4$.  Whether it's true that the Hermite constant in any rank $k$ is realized by a positive definite even lattice of minimum~$4$ is an open question.
\end{remark}

\section{Refined Artin--Tate conditions}
\label{sec:computation}

In this section, we use the sphere packing bounds in Section~\ref{sec:sphere}, together with the quadratic and cyclotomic Artin--Tate formulas, to establish new nontrivial conditions on the possible Weil polynomials of K3 surfaces over a finite field.  We call these \emph{refined Artin--Tate conditions}.

\begin{theorem}
\label{thm:refined_Artin-Tate_condition}
Let $X$ be a K3 surface over $\FF_q$ with Weil polynomial $P_X(t)$.  Fix $m \geq 2$ and write
\[
P_X(t) = (1-t)^\rho \, C_m(t) \, P_m(t)
\]
where $C_m(t) = \prod_{d \mid m,\, d>1} \Phi_d(t)^{a_d}$ for $a_d \geq 0$ and $P_m(\zeta_d) \neq 0$ for all $d \mid m$. Then we have
\[
|\NS_m(X)| = \frac{|A_m(X)|}{|\Br_m(X)|} \prod_{d \mid m,\, d>1} N_{\QQ(\zeta_d)/\QQ}(P_m(\zeta_d)).
\]
Concretely, we have that
\[
\Res\left(\frac{t^m-1}{t-1},P_m(t)\right) = \frac{|\Br_m(X)|\, |\NS_m(X)|}{|A_m(X)|}
\]
and when $\rho < \rho_m$, this yields the following refined Artin--Tate condition 
\[
\gcd\bigl(m^\rho, C_m(1)\bigr) \, \cdot \, \Res\left(\frac{t^m-1}{t-1},P_m(t)\right) \geq \delta_{\rho_m-\rho}
\]
where $\Res$ is the classical resultant and $\delta_r$ is the minimum discriminant of an even positive definite rootless lattice of rank $r$.
\end{theorem}

\begin{proof}
Applying M\"obius inversion to the first version of Theorem~\ref{thm:cyclotomic} yields
\[
\prod_{d \mid m} N_{\QQ(\zeta_d)/\QQ}(P_X^*(\zeta_d)) = m^{\rho_m}\, |\Br(X_m)| \, |\NS(X_m)|.
\]
Now write $G_m(t) = (1-t)^\rho \, C_m(t)$, and notice that
\[
\prod_{d \mid m} N_{\QQ(\zeta_d)/\QQ}(G_d^*(t)) = \prod_{\zeta^m=1} G_d^*(\zeta) = m^{\rho_m}
\]
since all roots of $G_m$ are $m$th roots of unity.  Since $P_X^*(\zeta) = G_m^*(\zeta)\, P_m(\zeta)$ for any $m$th root of unity $\zeta$, we obtain
\[
\prod_{d \mid m} N_{\QQ(\zeta_d)/\QQ}(P_m(\zeta_d)) = |\Br(X_m)| \, |\NS(X_m)|
\]
which can be viewed as a version of the Artin--Tate formula for $P_m(t)$.  The classical Artin--Tate formula gives $P_X^*(1) = |\Br(X)| \, |\NS(X)| = C_m(1)\, P_m(1)$.  Hence dividing both sides of the above by $P_m(1)$ yields
\[
\prod_{d \mid m,\, d>1} N_{\QQ(\zeta_d)/\QQ}(P_m(\zeta_d)) = \frac{C_m(1)\, |\Br(X_m)| \, |\NS(X_m)|}{|\Br(X)| \, |\NS(X)|} = \frac{|\Br_m(X)|}{|A_m(X)|}|\NS_m(X)|
\]
using the same identities as in the proof of Theorem~\ref{thm:cyclotomic}.  This is the desired formula.

For the concrete statement, we first remark that classical formulas for the resultant, cf.~\cite{apostol}, yield
\[
\Res(t^{m-1} + \dotsm + 1,P_m(t)) = \prod_{\zeta^m=1,\, \zeta \neq 1} P_m(\zeta_d) = \prod_{d \mid m,\, d>1} N_{\QQ(\zeta_d)/\QQ}(P_m(\zeta_d)). 
\]
Second, Lemma~\ref{lem:m} implies that $|A_m(X)|$ divides $b_m(X) = C_m(1)$ and Lemma~\ref{lem:invar_coinvar_basics}(e) implies that $|A_m(X)|$ divides $m^\rho$.  Third, $|\Br_m(X)|$ is an integer.  Thus we have that
\[
\gcd(m^\rho, C_m(1)) \Res(t^{m-1} + \dotsm + 1,P_m(t)) \geq |\NS_m(X)| \geq \delta_{\rho_m-\rho}
\]
since, by Lemma~\ref{lem:rootless}, $\NS_m(X)$ is (up to sign) a positive definite even rootless lattice, which has rank $\rho_m-\rho$, hence by Section~\ref{sec:sphere}, its discriminant is bounded by $\delta_{\rho_m-\rho}$.
\end{proof}

The introduction of the resultant in the statement of Theorem~\ref{thm:refined_Artin-Tate_condition} provides a more computationally practical way of computing the $\prod_{d \mid m,\, d>1} N_{\QQ(\zeta_d)/\QQ}(P_m(\zeta_d))$ term.  From McMullen's work on dynamics of automorphisms of K3 surfaces~\cite{mcmullen:dynamics,mcmullen:glue,mcmullen:entropy}, see also Brandhorst~\cite{brandhorst:nonsymplectic}, such resultants naturally arise as they directly control primes dividing the gluing groups of primitive extensions.  Hence the concrete formula in Theorem~\ref{thm:refined_Artin-Tate_condition} could likely be achieved up to taking squarefree parts using an adaptation of the techniques in \cite[Section~5]{brandhorst:nonsymplectic}.

As we shall see, these refined Artin--Tate conditions can be effectively used to rule out whether a Weil polynomial $P(t)$ of degree 22 can be realized by a K3 surface over $\FF_q$.  To sharpen the refined Artin--Tate conditions, we can apply them to a Weil polynomial $P(t) = q\prod_i (1-\alpha_i t)$ and to all its formal base changes $P^{(l)}(t) = q^l\prod_i (1-\alpha_i^l t)$.  In fact, there is only a finite number of possible base changes to check before the entire geometric Picard group is defined over $\FF_{q^l}$, and the conditions (which require further Picard rank jumping after an extension) fail to yield any further nontrivial restrictions.  Similarly, there is only a finite number of refined Artin--Tate conditions (a finite number of $m$) to check, governed by the possible factorizations of the algebraic part of the Weil polynomial into cyclotomic polynomials.

In the case $\rho = \rho_m$, the refined Artin--Tate condition doesn't apply.  However,  Corollary~\ref{cor:EJ_higher_m} gives a potential condition, namely that $N_{\QQ(\zeta_m)/\QQ}(P_1^*(\zeta_m))$ is a square.  Unfortunately, in view of Remark~\ref{rem:EJ_higher}, a result of Elsenhans and Jahnel~\cite[Corollary~18]{EJ} implies that this condition is automatically satisfied by any Weil polynomial when $m$ is odd, and is satisfied for all $m$ as soon as it is satisfied for $m=2$, reducing us to Elsenhans and Jahnel's original condition. 

\bigskip

We now review how the refined Artin--Tate conditions sort out the list of candidate Weil polynomials for K3 surfaces over $\FF_2$ computed by Kedlaya and Sutherland~\cite{KS}.

Let $\mathcal K$ be the set of $2{,}258{,}428$ candidate Weil polynomials of degree $22$ that satisfy only the point count nonnegativity and extension growth properties for K3 surfaces  over $\FF_2$, cf.\ \cite[Theorem~3.4]{AP}. We consider the subsets $E_2, R_2, R_{>2}$ of $\mathcal K$ consisting of Weil polynomials that fail the original Elsenhans and Jahnel quadratic Artin--Tate condition (when $\rho=\rho_2$), fail the quadratic refined Artin--Tate condition (when $m=2$), and fail the refined Artin--Tate condition for some $m>2$, respectively.  As discussed above, in computing these subsets, we also take into account evaluating each condition after performing formal base changes to $\FF_{q^l}$.  In particular, $\mathcal K \smallsetminus E_2$ is the set of $1{,}672{,}565$ candidate Weil polynomials of K3 surfaces over $\FF_2$ in \cite[Computation~3c]{KS}.  See Table~\ref{tab:at-family-exclusions} for the statistical summary.

\bigskip

\begin{table}[htbp]
\centering
\begin{tabular}{@{}lrrr@{}}
\toprule
Subset  & Total in $\mathcal K$ & New beyond $E_2$ & With $E_2$\\
$S \subset \mathcal K$  & $|S|$ & $|S\setminus E_2|$ & $|S\cup E_2|$ \\
\midrule
Original E-J quadratic $(E_2)$ & $585{,}863$ & $0$ & $585{,}863$ \\
Refined quadratic $(R_2)$ & $770{,}118$ & $628{,}290$ & $1{,}214{,}153$ \\
Refined cyclotomic $(R_{>2})$ & $870{,}749$ & $658{,}298$ & $1{,}244{,}161$ \\
\midrule
Total & $1{,}464{,}443$ & $878{,}580$ & $1{,}464{,}443$ \\
\bottomrule\smallskip
\end{tabular}
\caption{Exclusion counts for the three types of Artin--Tate conditions in $\mathcal K$.}
\label{tab:at-family-exclusions}
\end{table}

\newpage

In Kedlaya and Sutherland's candidate list $\mathcal K \smallsetminus E_2$, the condition $R_{>2}$ excludes $250{,}290$ candidates that survive every quadratic condition. Conversely, in $\mathcal K \smallsetminus E_2$, the quadratic conditions $E_2, R_2$ exclude $220{,}282$ candidates that survive every higher degree cyclotomic condition. Together, the new conditions exclude $878{,}580$ of the $1{,}672{,}565$ Kedlaya--Sutherland candidates ($52.53\%$), leaving a candidate list of $793{,}985$ Weil polynomials.

\raggedbottom
\bibliographystyle{amsplain}
\bibliography{references}
\end{document}